\documentclass[11pt,letterpaper]{amsart}
\usepackage[left=20mm,right=20mm,top=20mm,bottom=20mm]{geometry}
\pdfoutput=1
\usepackage{graphicx}
\usepackage{amsfonts,amsthm,amsmath,amssymb,latexsym, amscd, euscript}
\usepackage{graphicx,color}
\usepackage[all]{xy}
\usepackage{epsfig}
\usepackage{mathrsfs}
\DeclareMathOperator{\sech}{sech}
\DeclareMathOperator{\csch}{csch}

\DeclareUnicodeCharacter{2212}{-}

\usepackage{amssymb}
\usepackage{epstopdf}
\theoremstyle{plain}
\newtheorem{theorem}{Theorem}[section]

\newtheorem{definition}[theorem]{Definition}

\begin{document}

\theoremstyle{definition} 

\newtheorem*{notation}{Notation}  

\theoremstyle{plain}      

\def\H{{\mathbb H}}
\def\F{{\mathcal F}}
\def\R{{\mathbb R}}
\def\Q{\hat{\mathbb Q}}
\def\Z{{\mathbb Z}}
\def\E{{\mathcal E}}
\def\N{{\mathbb N}}
\def\X{{\mathcal X}}
\def\Y{{\mathcal Y}}
\def\C{{\mathbb C}}
\def\D{{\mathbb D}}
\def\G{{\mathcal G}}
\def\T{{\mathcal T}}

\title{Non-ergodicity of the geodesic flow on Cantor tree surfaces}


\date{}
\author{Michael Pandazis}

\address[Michael Pandazis]{PhD Program in Mathematics, The Graduate Center, CUNY \\ 365 Fifth Ave., N.Y., N.Y., 10016, USA.}
\email{mpandazis@gradcenter.cuny.edu}

\maketitle
\vspace{-.5in}
\begin{abstract}
    A Riemann surface equipped with its conformal hyperbolic metric is parabolic if and only if the geodesic flow on its unit tangent bundle is ergodic. Let $X$ be a Cantor tree or a blooming Cantor tree Riemann surface. Fix a geodesic pants decomposition of $X$ and call the boundary geodesics in the decomposition {\it cuffs}. Basmajian, Hakobyan, and \v Sari\' c proved that if the lengths of cuffs are rapidly converging to zero, then $X$ is parabolic. More recently, \v Sari\' c proved a slightly slower convergence of lengths of cuffs to zero implies $X$ is not parabolic. In the paper, we interpolate between the two rates of convergence of the cuffs to zero and find that these surfaces are not parabolic, thus completing the picture.
\end{abstract}

\vspace{-.40cm}
\section{Introduction}

    A Riemann surface $X$ is parabolic, denoted by $X \in O_G$, if it does not admit a Green's function-i.e., a harmonic function $u:X\to\mathbb{R}^{+}$ with a logarithmic singularity at a single point of $X$ whose values limit to zero at the ideal boundary (Ahlfors-Sario \cite{Ahlfors}). It is known that $X\in O_G$ if and only if the geodesic flow (for the conformal hyperbolic metric) on the unit tangent bundle of $X$ is ergodic if and only if the Poincar\'e series for the covering Fuchsian group diverges if and only if the Brownian motion on $X$ is recurrent (see Nicholls \cite{Nicholls}, Sullivan \cite{Sullivan}, Tsuji \cite{Tsuji}, Basmajian-Hakobyan-\v{S}ari\'{c} \cite{BHS}).

When $X$ is of finite type, then $X\in O_G$ if and only if $X$ has finite area. A Riemann surface $X$ is said to be {\it infinite} if its fundamental group cannot be finitely generated. An infinite Riemann surface is determined by a fixed geodesic pants decomposition and the Fenchel-Nielsen parameters associated to the pants decomposition (Basmajian-\v{S}ari\'{c} \cite{BasmajianSaric}). As in \cite{BHS}, we consider the question of deciding when $X\in O_G$ based on its Fenchel-Nielsen parameters.
\vspace{.35in}
\begin{figure}[ht]
\begin{picture}(100,100)
\put(-25,0){\includegraphics[width = 2in]{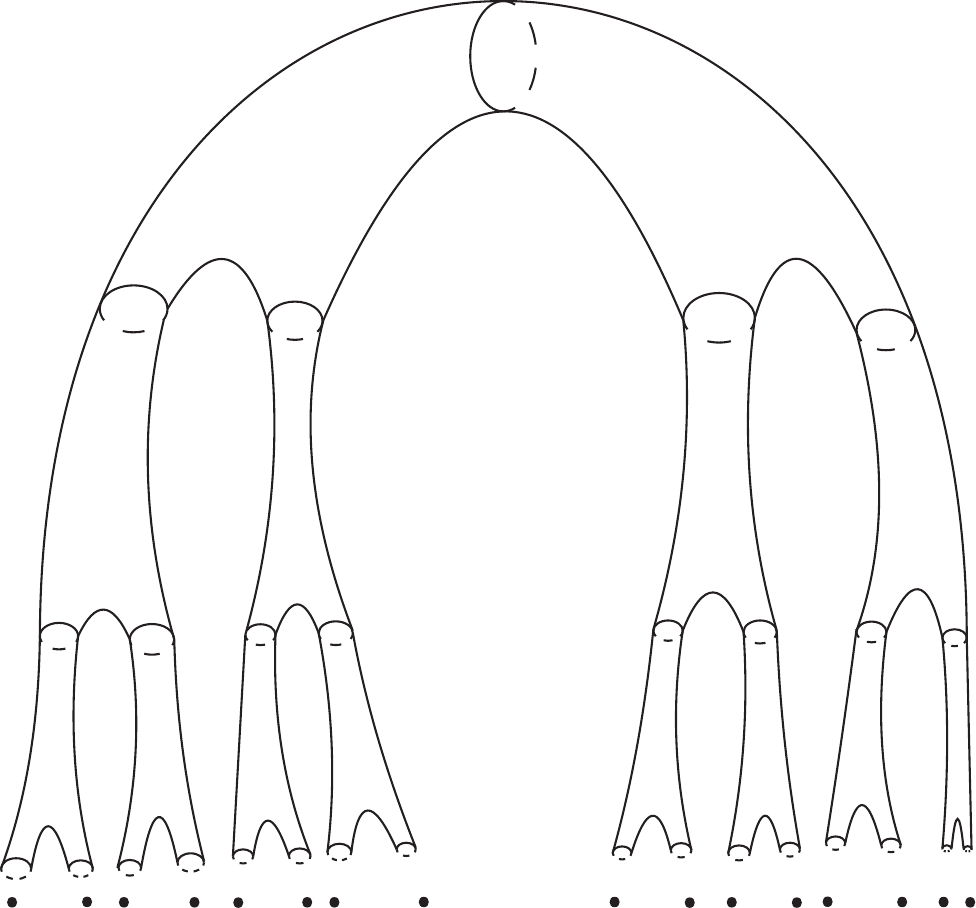}}
\put(112,94){$X_1$}
\put(123,45){$X_2$}
\put(124,12){$X_3$}
\end{picture}
\vspace{-.5cm}
\caption{The Cantor tree surface with a geodesic pants decomposition.}
\label{fig:The Cantor tree surface}
\end{figure}
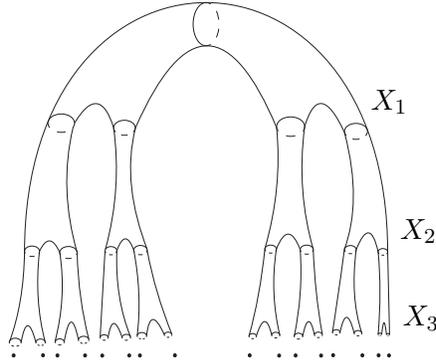

A Cantor tree Riemann surface $X_C$ is conformal to the complement of a Cantor set in the Riemann sphere. Equivalently, $X_C$ is constructed by isometrically gluing countably many geodesic pairs of pants along their boundary geodesics (called {\it cuffs}) to form the ``shape'' of the dyadic tree (see Figure \ref{fig:The Cantor tree surface}). In addition to the lengths of its cuffs, the Cantor tree Riemann surface $X_C$ is determined by the twists along the cuffs.

The cuffs of $X_C$ are grouped in the levels based on the level in the dyadic tree. At level zero, we have a single cuff, which is at the top of $X_C$ in Figure 1. At level one, we have four cuffs, and at level $n\geq 1$, we have $2^{n+1}$ cuffs. Denote by $\{\alpha_n^j\}_{j=1}^{2^{n+1}}$ the level $n$ cuffs from left to right in Figure \ref{fig:The Cantor tree surface}.

McMullen \cite{McMullen} proved if there is a $C>0$ such that $1/C\leq \ell( \alpha_n^j)\leq C$ then $X\notin O_G$. This is a consequence of the fact that the Brownian motion has many directions to escape to infinity when the ideal boundary is large, in our case the Cantor set, and the cuffs (the openings) are not short.
In the case when the cuffs are short Basmajian, Hakobyan and \v Sari\' c \cite{BHS} proved $X_C\in O_G$ if there is $C>0$ such that $$\ell (\alpha_n^j)\leq C\frac{n}{2^n},$$ where $\ell (\cdot )$ is the hyperbolic length in $X$. 

More recently, \v{S}ari\'{c} \cite[Theorem 8.3]{Saric} proved that if 
$$
\ell( \alpha_n^j)=\frac{n^r}{2^n}
$$
for $r>2$, and for all $n\geq 1$ and $j=1,2,\ldots ,2^{n+1}$ then $X\notin O_G$. Thus, the Brownian motion escapes to infinity even when the cuffs are short in this controlled fashion.

The remaining case to consider is whether $X_C$ is parabolic or not for $1<r\leq 2$. We show the following.

\begin{theorem}
\label{thm:Cantor_trees_non_parabolic}
Let $X_C$ be the Cantor tree surface as depicted in Figure \ref{fig:The Cantor tree surface} and $\{\alpha_n^j\}_{j=1}^{2^{n+1}}$ the cuffs at the level $n$. The cuff lengths are decreasing along each end. Then $X_C\notin O_G$ if there is an $r>1$ such that
$$
C_1\frac{n^r}{2^n}\leq\ell (\alpha_n^j )\leq \frac{C_2}{n^2}
$$
for some universal constants $C_1,C_2>0$.
\end{theorem}
$ \newline$
\vspace{.50in}
\begin{figure}[ht]
\begin{picture}(100,100)
\put(-55,0){\includegraphics[width = 2.8in, height=2.4in]{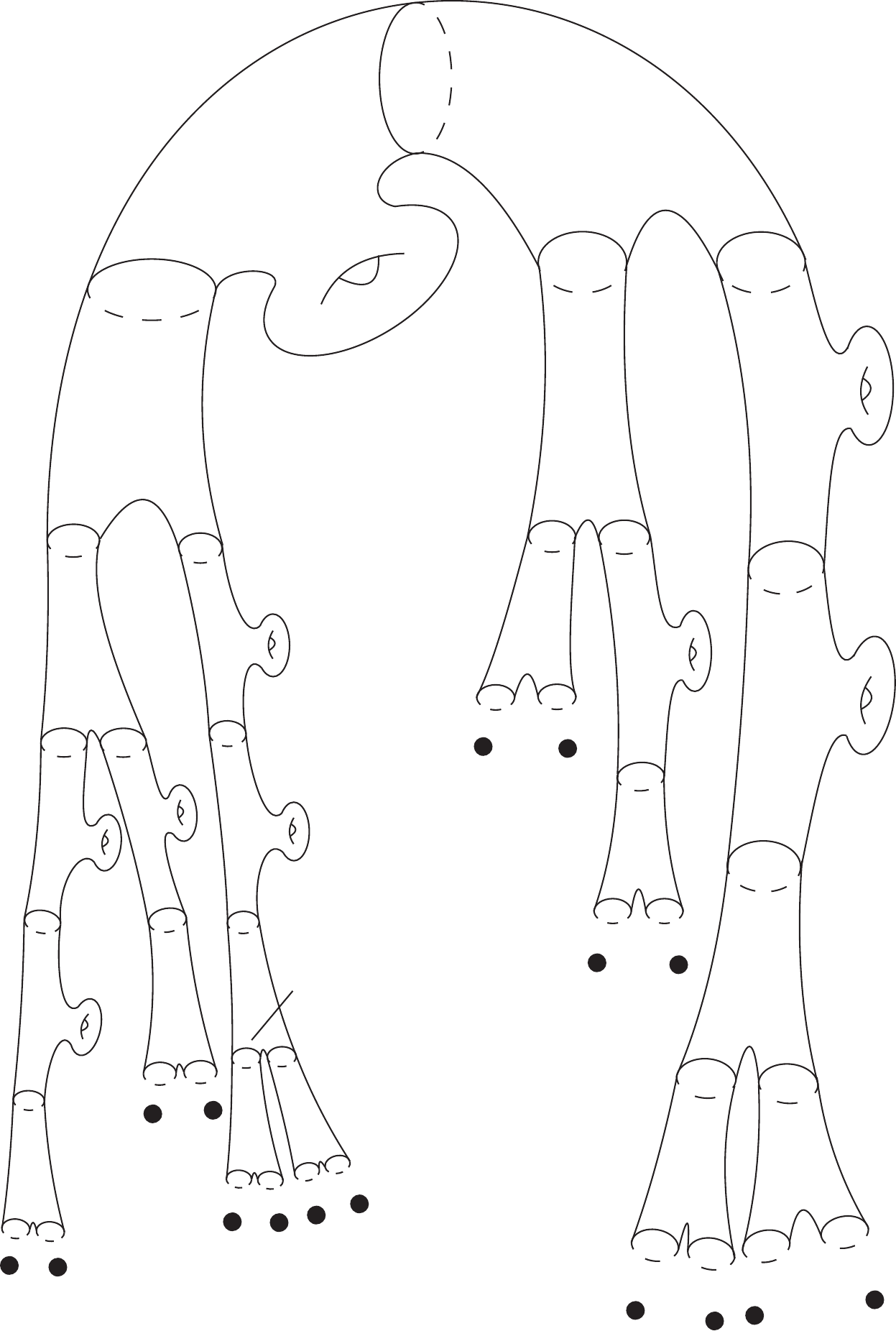}}
\put(33,161){$\alpha_0$}
\put(-51,136){$\alpha_0$}
\put(-58,100){$\alpha_1^1$}
\put(5,52){$\alpha_1^2$}
\put(50,134){$\alpha_1^3$}
\put(129,58){$\alpha_1^4$}
\put(-67,28){$\alpha_2^1$}
\put(-36,50){$\alpha_2^2$}
\put(12,43){$\alpha_2^3$}
\put(16,32){$\alpha_2^4$}
\put(49,102){$\alpha_2^5$}
\put(96,70){$\alpha_2^6$}
\put(83,31){$\alpha_2^7$}
\put(132,30){$\alpha_2^8$}
\end{picture}
\vspace{-.5cm}
\caption{The blooming Cantor tree surface with a geodesic pants decomposition.}
\label{fig:The blooming Cantor tree surface}
\end{figure}
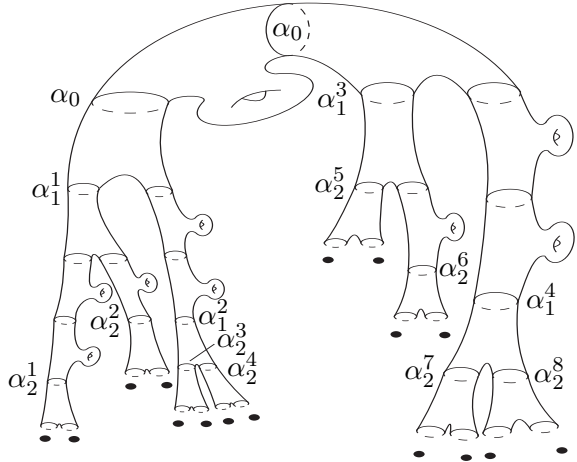
\\
\indent Even for $r>2$, the scope of our theorem is slightly more general than \cite[Theorem 8.3]{Saric} because we allow the lengths of the cuffs to vary with the given lower bound. We also extend our result to surfaces with infinite genus and a Cantor set of ends, called the {\it blooming Cantor tree surfaces} $\tilde{X}_C$ (see Figure \ref{fig:The blooming Cantor tree surface}). To construct the blooming Cantor tree from the Cantor tree, attach a geodesic surface of genus at most $C$ and two boundaries, or do not, to each level $n$ boundary $\alpha_n^j$ (see that $\alpha_0^1=\alpha_0^2=\alpha_0$) and redefine $\alpha_n^j$ to be the boundary of the attached surface further away from $\alpha_0$ for $n \ge 0$ and for $1 \le j \le 2^{n+1}$ and for a universal constant $C>0$. We can add $2^{n+1}$ surfaces at the level $n$ of genus at most $C$ for $n \ge 0$. Assume the lengths of cuffs along each attached geodesic surface are decreasing.


\begin{theorem}
\label{thm:blooming cantor-tree non parabolicity}
Let $\tilde{X}_C$ be the blooming Cantor tree surface and $\{\alpha_n^j\}_{j=1}^{2^{n+1}}$ the cuffs at the level $n$ as depicted in Figure \ref{fig:The blooming Cantor tree surface}. The cuff lengths are decreasing along each end and each geodesic subsurface between level $n$ and level $n+1$ boundaries of $\tilde{X}_C$ has genus bounded above by $C>0$. Then $\tilde{X}_C\notin O_G$ if there is an $r > 1$ such that
$$
C_1\frac{n^r}{2^n}\leq\ell (\alpha_n^j )\leq \frac{C_2}{n^2}
$$
for some universal constants $C_1,C_2 >0$.
\end{theorem}

\begin{center}
    ACKNOWLEDGEMENTS
\end{center}

\indent The author would like to thank Dragomir \v Sari\' c for asking a question that led to the paper and for his advice and guidance. Additionally, a special thanks to Ara Basmajian, Hrant Hakobyan, and Dragomir \v Sari\' c, whose work on the type problem motivated the proof of non-ergodicity of the geodesic flow on Cantor tree surfaces.

\section{Partial measured foliations and a sufficient condition for a surface to be non-parabolic}
\indent Let $X=\mathbb{H}/\Gamma$ be an infinite Riemann surface, where $\mathbb{H}$ is the hyperbolic plane and $\Gamma$ is a Fuchsian covering group. For our purposes, a special case of the definition in \cite{Saric} with $E_i=U_i$ is enough.
\begin{definition}[\cite{Saric}]
    A {\it partial measured foliation} $\mathscr{F}$ on $X$ is an assignment of a collection of sets $\{U_i\}_i$ of $X$ (which do not have to cover the entire surface $X$) and continuously differentiable (with surjective tangent map) real-valued functions
$$
v_i: U_i \rightarrow \mathbb{R}.
$$
The sets $U_i$ are closed Jordan domains with piecewise differentiable boundaries. The pre-image $v_i^{-1} (c)$ for $c \in \mathbb{R}$ is a connected differentiable arc with endpoints on $\partial U_i$, and
\begin{equation}
\label{eq:v_i-on-intersection-of-open-sets}
v_i = \pm v_j + const
\end{equation}
on $U_i \cap U_j$. The collection of sets $\{U_i\}_i$ is {\it locally finite} in $X$. 
\end{definition}


 Let $\mathscr{F}$ be a fixed partial measured foliation on the surface $X$. A curve in $X$ is said to be a {\it horizontal arc} if it is expressible as a finite or infinite connected union of curves defined by $v_i^{-1} (c_i)$ for some collection of real numbers $c_i$. When a curve in $X$ is a maximal horizontal arc, it is called a {\it horizontal trajectory} of $\mathscr{F}$. A partial measured foliation is {\it proper} if each end of the lift to the universal cover $\mathbb{H}$ of every horizontal trajectory approaches a distinct point on the ideal boundary of $\mathbb{H}$.

\indent The {\it Dirichlet integral} (see \cite{Ahlfors}) of a continuously differentiable function $v_i:U_i\to\mathbb{R}$ is
\begin{equation}
\label{eq:dirichlet integral}
\int_{U_i} [(\frac{\partial v_i}{\partial x})^2+(\frac{\partial v_i}{\partial y})^2]dxdy.
\end{equation}
\indent The Dirichlet integral $D_X(\mathscr{F})$ of $\mathscr{F}$ over $X$ is
$$
D_X(\mathscr{F}) = \sum_{i}\int_{U_i} [(\frac{\partial v_i}{\partial x})^2+(\frac{\partial v_i}{\partial y})^2]dxdy,
$$
when the $U_i$'s are non-overlapping sets up to a set of measure zero.
A proper partial measured foliation $\mathscr{F}$ on $X$ is {\it integrable} if $D_X (\mathscr{F}) < \infty$.

 From \cite[Theorem 3.3]{Saric} and \cite[Theorem 4.1]{Saric} 
 it immediately follows
\begin{theorem}
\label{thm:section_2_theorem}
If there is a non-trivial integrable partial measured foliation of a Riemann surface $X$ with leaves that escapes every compact subset of $X$ at both ends, then $X$ is not parabolic.
\end{theorem}

\section{Proof of Theorems 1.1 and 1.2}
$ \newline$
\vspace{-.1in}
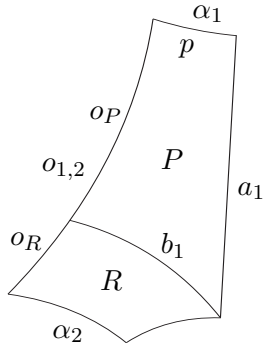
\begin{figure}[ht]
\begin{picture}(100,100)
\put(-10,0){\includegraphics[width = 1.2in, height=1.7in]{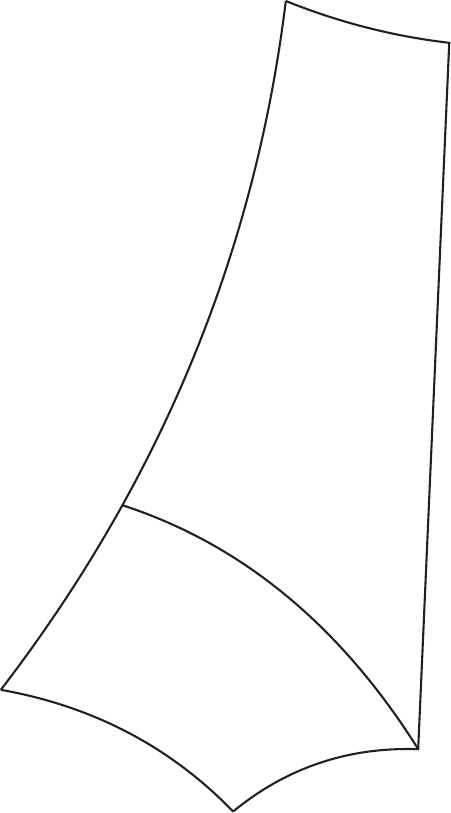}}
\put(60,123){$\alpha_1$}
\put(7,2){$\alpha_2$}
\put(48,66){$P$}
\put(25,20){$R$}
\put(21,85){$o_P$}
\put(-9,37){$o_R$}
\put(3,65){$o_{1,2}$}
\put(48,34){$b_1$}
\put(77,56){$a_1$}
\put(55,110){$p$}
\end{picture}
\vspace{-.5cm}
\caption{Horizontal foliation through $P \cup R$.}
\label{fig:Horizontal Foliation through Omega one and two}
\end{figure}
\vspace{-.06cm}

Define {\it a to be asymptotic to b}, denoted by $a \asymp b$, to mean that there is a $k>0$ such that $\frac{1}{k} \le \frac{a}{b} \le k$. Define {\it a to be asymptotically less than b}, notated by $a \lesssim b$, to mean that there is a $k>0$ such that $\frac{a}{b} \le k$.

It is enough to construct an integrable partial measured foliation on $X_C$ when all twists are zero because varying the twists by a bounded amount is a quasiconformal deformation \cite{ALPS} and parabolicity is a quasiconformal invariant \cite{Ahlfors}. Each geodesic pair of pants is divided into two right-angled hexagons by three orthogeodesic arcs between the pairs of cuffs. Since all twists are zero, the union of the orthogeodesic arcs forms a family of bi-infinite geodesics that separates $X_C$ into two symmetric halves permuted by an orientation-reversing isometry (see Figure 1).

Consider a pair of pants $\Pi$ from the decomposition with boundaries $\alpha_1$, $\alpha_2$, and $\alpha_3$. Let $o_{i,j}$ be the orthogeodesic arc between $\alpha_i$ and $\alpha_j$, for $i,j \in \{1,2,3\}$ such that $i \neq j$. The union $o_{1,2} \cup o_{1,3} \cup o_{2,3}$ separates $\Pi$ into front and back hexagons $H_1$ and $H_2$ with geodesic boundaries. Let $a_1$ be the orthogeodesic from $\alpha_1$ to $o_{2,3}$ that separates $H_1$ into two right-angled pentagons and divides $\alpha_1$ in $H_1$ from left to right into the sub-arcs $p$ and $q$. Call $P_p$ the pentagon containing $p$ (see Figure \ref{fig:Horizontal Foliation through Omega one and two}), and $P_q$ the pentagon containing $q$. Since the lengths of $\alpha_2$ and $\alpha_3$ are not necessarily the same, $p$ is not necessarily equal to $q$.

The orthogeodesic $b_1$ from $a_1 \cap o_{2,3}$ to $o_{1,2}$ divides $P_p$ into quadrilaterals $P$ and $R$ adjacent to $\alpha_1$ and $\alpha_2$. Let $o_P$ be the sub-arc of $o_{1,2}$ in $P$ and let $o_R$ be the sub-arc of $o_{1,2}$ in $R$ (see Figure \ref{fig:Horizontal Foliation through Omega one and two}). Lift $P$ isometrically to $\mathbb{H}$ as follows. Lift the geodesic arc $o_P$, starting at $\alpha_1$, to the y-axis from $i$ to $e^{\ell(o_P)}i$ and call it $\tilde{o}_P$. Each point $w$ in $P$ belongs to a hyperbolic geodesic arc $\gamma_{w_0}$ orthogonal to $o_P$ with foot $w_0$ on $o_P$. Map $\gamma_{w_0}$ to the geodesic arc $\tilde{\gamma}_{w_0}$ orthogonal to the y-axis in $\mathbb{H}$, to its right, and whose foot on the y-axis is the lift $\tilde{w}_0$ of $w_0$ for each $w_0$ in $o_P$. That defines an isometric lift $\tilde{P}$ of $P$ to $\mathbb{H}$. Denote the lift of $a_1$ by $\tilde{a}_1$.
$ \newline$
\vspace{-.2in}
\begin{figure}[ht]
\begin{picture}(100,100)
\put(-5,0){\includegraphics[width=1.5in]{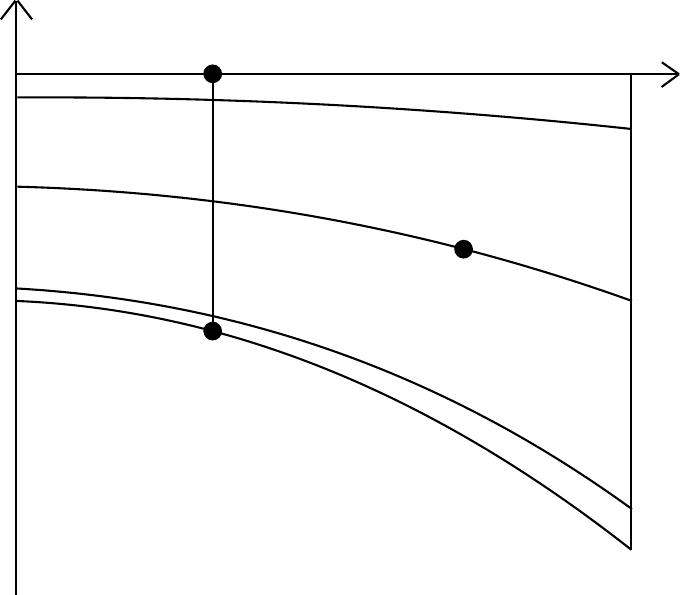}}
\put(10,87){$x$}
\put(45,64){$z = x + iy$}
\put(32,48){$d_{h}^{P}(x)$}
\put(-5,100){$y$}
\put(-34,80){$v = 0$}
\put(106,81){$x$}
\put(-34,42){$v = 1$}
\end{picture}
\vspace{-.5cm}
\caption{The image $f(\tilde{P})$ of a lift of Lambert quadrilateral $P$ in $X$ by a diffeomorphism $f$.}
\label{fig:image of Lambert quadrilateral to a subset of complex plane}
\end{figure}
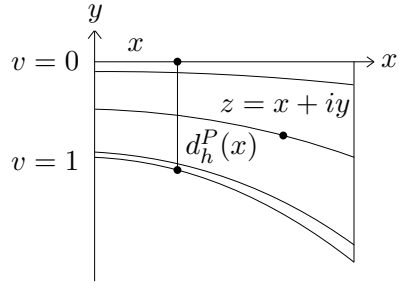
\\
\indent Map $\tilde{o}_P$ to $[0,\ell(o_P)]$ on $\mathbb{R}$ and each $\tilde{\gamma}_{w_0}$ to a Euclidean segment orthogonal to $[0,\ell(o_P)]$ and below the x-axis by an isometry. That defines $f: \tilde{P}\to \mathbb{C}$ (see Figure \ref{fig:image of Lambert quadrilateral to a subset of complex plane}). The length of $f(\tilde{\gamma}_{w_0})$, where $f(\tilde{w}_0) = x$, by a formula for Lambert quadrilaterals from \cite[Theorem 2.3.1(iv)]{Buser}  is
\begin{equation}
\label{eq:d_h^1}
d_h^P(x)=\tanh^{-1}(\cosh x\tanh \ell(p)).
\end{equation}

Define a real-valued, continuously differentiable function $v_P$ with $d_h^P (x)$ from (\ref{eq:d_h^1}) to be
$$
v_P (x+iy) = \frac{y}{-d_h^P (x)} = \frac{-y}{
\tanh^{-1}(\cosh x\tanh \ell(p))
}
$$
for $z=x + iy$ in $f(\tilde{P})$.
The function $v = v_P$ defines a horizontal foliation in $f(\tilde{P})$ with leaves defined by $v_P^{-1}(c)$ for $0 \le c \le 1$ (see Figure \ref{fig:image of Lambert quadrilateral to a subset of complex plane}). We obtain an upper estimate of the integrand of the Dirichlet integral over $f(\tilde{P})$ and then use it to estimate the integral for sufficiently small lengths of $\alpha_1$, $\alpha_2$, and $\alpha_3$. By Lemma A.4(2), $d/dx(\tanh^{-1}x)=1/(1- x^2)$, $\tanh^{-1}x\geq x$ for $x\geq 0$ and $\ell (p)\to 0$,
$$
(\frac{\partial v_P}{\partial x})^2 = y^2[\frac{1}{\tanh^{-1}{(\tanh{\ell(p)}\cosh{x})}}]^4[\frac{1}{1 - [\tanh{\ell(p)}\cosh{x}]^2}]^2 \tanh^2{\ell(p)}\sinh^2{x}
$$
$$
\lesssim \frac{y^2}{\tanh^2{\ell(p)}}\tanh^2{x}\sech^2{x} \lesssim \frac{y^2}{\ell(p)^2}\tanh^2{x}\sech^2{x},
$$
and
$$
(\frac{\partial v_P}{\partial y})^2 = [\frac{1}{\tanh^{-1}(\tanh{\ell(p)}\cosh{x})}]^2 \lesssim \frac{1}{\ell(p)^2}\sech^2{x}.
$$
By \cite[Theorem 2.3.4(ii)]{Buser} and $\sinh{(\cosh^{-1}{x})} < x$ for all $x > 1$,
\begin{equation}
\label{eq:1}
\ell(p)\sinh{\ell(o_{1,2})} = \ell(p)\sinh{(\cosh^{-1}{(\frac{1}{\tanh{\frac{\ell(\alpha_2)}{2}}\tanh{\ell(p)}})})} \lesssim \frac{1}{\ell(\alpha_2)}.
\end{equation}
We integrate, use equations (\ref{eq:d_h^1}) and (\ref{eq:1}), and  Lemma A.4(3) to get
$$
\iint_{f(\tilde{P})}(\frac{\partial v_P}{\partial x})^2 \lesssim \frac{1}{\ell(p)^2}\int_{0}^{\ell(o_P)} \int_{0}^{d_h^{P}(x)} y^2\tanh^2{x}\sech^2{x}dydx
$$
$$
\lesssim \ell(p)\int_{0}^{\ell(o_P)} \sinh{x}\tanh{x}dx < \ell(p)\int_{0}^{\ell(o_P)} \cosh{x}dx < \ell(p)\sinh{\ell(o_{1,2})} \lesssim \frac{1}{\ell(\alpha_2)}.
$$
In addition, by Lemma A.4(3) and $\int_0^{\infty}\sech x dx = \frac{\pi}{2}$ we get
$$
\iint_{f(\tilde{P})}(\frac{\partial v_P}{\partial y})^2 \lesssim \frac{1}{\ell(p)^2}\int_{0}^{\ell(o_P)} \int_{0}^{d_h^{P}(x)} \sech^2{x}dydx \lesssim \frac{1}{\ell(p)}\int_{0}^{\ell(o_P)}\sech{x}dx \lesssim \frac{1}{\ell(p)}.
$$
The above, together with Lemma A.4(4), gives
$$
\iint_{f(\tilde{P})} (\frac{\partial v_P}{\partial x})^2 + (\frac{\partial v_P}{\partial y})^2 \lesssim \frac{1}{\ell(\alpha_2)} + \frac{1}{\ell(p)} \lesssim \frac{1}{\ell(\alpha_1)} + \frac{1}{\ell(\alpha_2)}.
$$

The orthogeodesic from the point $a_1 \cap o_{2,3}$ to $o_{1,3}$ divides $P_q$ into quadrilaterals $Q$ and $S$ adjacent to $\alpha_1$ and $\alpha_3$. By the analogous notation and derivations, we obtain
\begin{equation}
\label{eq:integral-over-Omega_1-union-Omega_2}
\iint_{f(\tilde{P}) \cup f(\tilde{R})} (\frac{\partial v}{\partial x})^2 + (\frac{\partial v}{\partial y})^2 \lesssim \frac{1}{\ell(\alpha_1)} + \frac{1}{\ell(\alpha_2)} \ \mathrm{and} \ \iint_{f(\tilde{Q}) \cup f(\tilde{S})} (\frac{\partial v}{\partial x})^2 + (\frac{\partial v}{\partial y})^2 \lesssim \frac{1}{\ell(\alpha_1)} + \frac{1}{\ell(\alpha_3)}.
\end{equation}

Lemma B.1 enables us to estimate the Dirichlet integrals of the foliations of quadrilaterals $P$, $Q$, $R$, and $S$ in the front of $\Pi$ from above using the inequality \cite{LarsAhlfors}
\begin{equation}
\label{eq: dirichlet integral comparison}
\iint_{\tilde{\Omega}} [(\frac{\partial (v \circ f)}{\partial \xi})^2 + (\frac{\partial (v \circ f)}{\partial \eta})^2] d \xi d \eta \le k_0 \iint_{f(\tilde{\Omega})} [(\frac{\partial v}{\partial x})^2 + (\frac{\partial v}{\partial y})^2] dx dy,
\end{equation}
where $\tilde{\Omega}$ is $\tilde{P}$, $\tilde{Q}$, $\tilde{R}$, or $\tilde{S}$.

\vskip .2 cm

We define an integrable partial measured foliation $\mathscr{F}$ supported on the front of the Cantor tree surface $X_C$ by scaling the partial foliations defined by $v_P$, $v_Q$, $v_R$, and $v_S$ in every pair of pants in the decomposition of $X_C$ in order for the transverse measures on the common boundaries of any two quadrilaterals to agree. The transverse measures on the ``vertical'' boundaries of the quadrilaterals $P$, $Q$, $R$, and $S$ that are given by integrating the differentials $dv_P$, $dv_Q$, $dv_R$, and $dv_S$ are proportional to the hyperbolic lengths, and on each vertical boundary, the corresponding measure equals to $1$.



\vspace{1.7in}
\begin{figure}[ht]
\begin{picture}(100,100)
\put(-75,0)
{\includegraphics[width=3.5in]{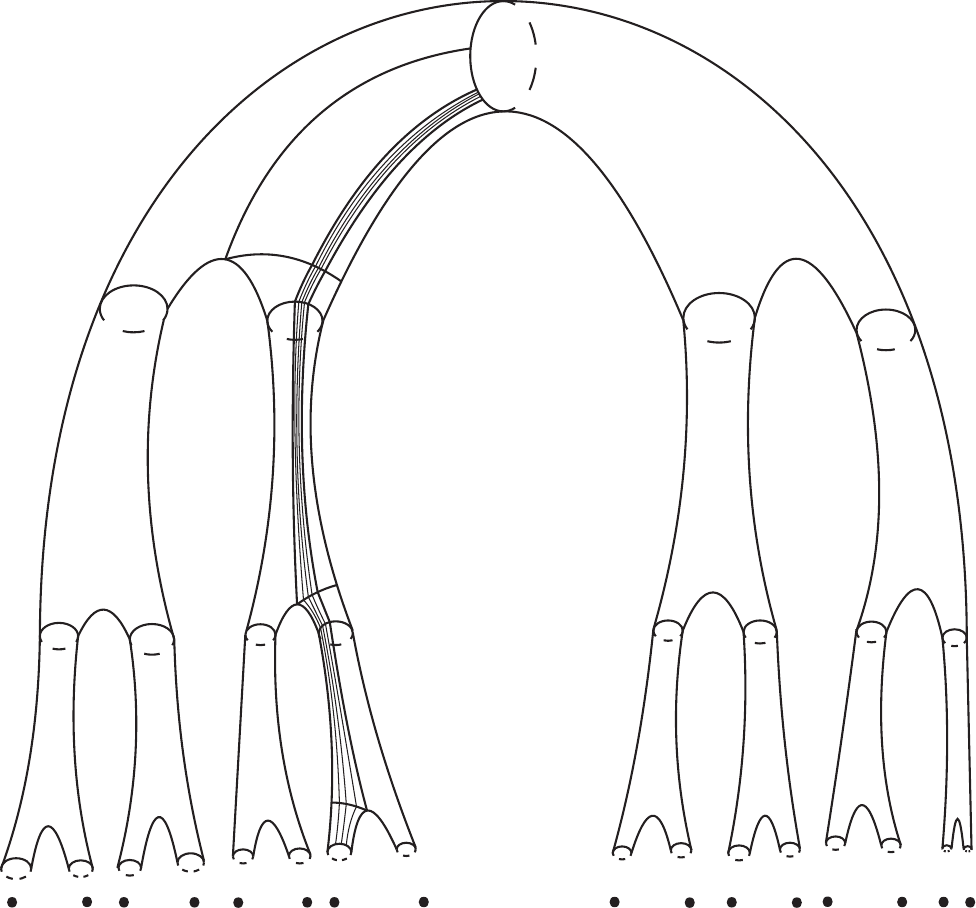}}
\put(49,214){$q_0$}
\put(13,149){$q_1^2$}
\put(-3,64){$p_2^4$}
\end{picture}
\vspace{-.5cm}
\caption{Illustration of how to choose $p$'s and $q$'s for the product of relative lengths.}
\label{fig:how_to_choose_p_q}
\end{figure}
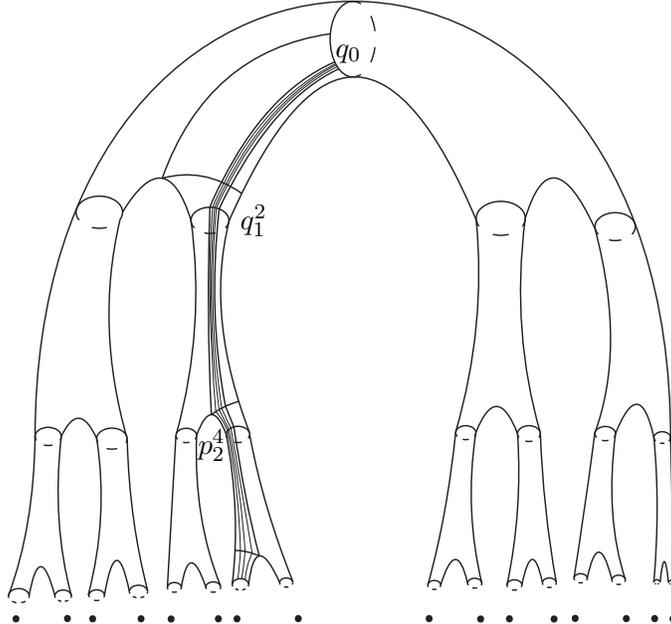

Consider any level $n+1$ pair of pants in $X$ with level $n$ boundary geodesic $\alpha_n^j$ and level $n+1$ boundary geodesics $\alpha_{n+1}^{2j-1}$ and $\alpha_{n+1}^{2j}$. Define the geodesic arcs $p_n^j$ and $q_n^j$ to be in the front hexagon of the pair of pants as the intersections $P \cap \alpha_n^j$ and $Q \cap \alpha_n^j$. We define the {\it relative lengths} of $p_n^j$ and $q_n^j$ by
$$
\ell_{n,j}^0=\frac{\ell (p_n^j)}{\ell (\alpha_n^j)/2}\ \ \mathrm{and}\ \ \ \ell_{n,j}^1=\frac{\ell (q_n^j)}{\ell (\alpha_n^j)/2}.
$$


Let $P^{n,j}$, $Q^{n,j}$, $R^{n,j}$, and $S^{n,j}$ be the Lambert quadrilaterals $P$, $Q$, $R$, and $S$ in the front of a level $n$ pair of pants $\Pi^{n,j}$ of $X$ with the front of a level $n-1$ cuff $\alpha_{n-1}^j$ as a boundary for $n \ge 1$ and for each $1 \le j \le 2^n$ (see that $\alpha_0^1 = \alpha_0^2 = \alpha_0$). Let $v_P^{n,j}$, $v_Q^{n,j}$, $v_R^{n,j}$, and $v_S^{n,j}$ be the partial measured foliations of $f(\tilde{P}^{n,j})$, $f(\tilde{Q}^{n,j})$, $f(\tilde{R}^{n,j})$, and $f(\tilde{S}^{n,j})$ for $n \ge 1$ and for each $1 \le j \le 2^{n}$ as the foliations $v_P$, $v_Q$, $v_R$, and $v_S$ for $f(\tilde{P})$, $f(\tilde{Q})$, $f(\tilde{R})$, and $f(\tilde{S})$.

Each pair of pants $\Pi^{n,j}$ starting from the top cuff $\alpha_0$ can be reached by a unique path of $n$ consecutive cuffs. In addition, at each cuff in the path, with the exception of the last cuff, we can choose either $p$ or $q$ depending on whether the next cuff is to the left or the right. The new function $w$ that defines the partial measured foliation $\mathscr{F}$ is obtained by multiplying the foliations of $f(\tilde{P}^{n,j})$ and $f(\tilde{R}^{n,j})$ with the product of the relative lengths of the corresponding choices of $p$'s and $q$'s on the path of cuffs times the relative length of $p_n^j$ (see Figure \ref{fig:how_to_choose_p_q}), and by multiplying the foliations of $f(\tilde{Q}^{n,j})$ and $f(\tilde{S}^{n,j})$ with the product of the relative lengths of the corresponding choices of $p$'s and $q$'s on the path of cuffs times the relative length of $q_n^j$. In this fashion, the transverse measures of the foliations of adjacent quadrilaterals on the common side of the quadrilaterals given by the function $w$ are equal, and $w$ defines a partial measured foliation supported on the front side of $X_C$.

The function $w$ induces the measure on $\alpha_0$ of mass $1$. 
The total mass of transverse measure on $\alpha_{n}^j$ is
\begin{equation}
\label{eq:transverse measure on alpha_0}
    \Pi_{k=0}^{n-1}\ell_{k,j}^{i_k},
\end{equation}
where $i_k\in\{ 0,1\}$ depending on the path of consecutive cuffs from $\alpha_0$ to $\alpha_n^j$. 

Let $T = \sum_{n=1}^{\infty}\frac{1}{n^2} < \infty$. From Appendix A.2,
$$
\frac{1}{2}e^{-\frac{C_2}{(k+1)^2}} \le \ell_{k,j}^{i_k} \le \frac{1}{2}e^{\frac{C_2}{(k+1)^2}}
$$
for $0 \le k \le n-1$. Use the above inequalities for $0 \le k \le n-1$ to obtain
$$
\frac{1}{2^n}e^{-C_2 T} \le \frac{1}{2^n}e^{-C_2\sum_{k=1}^{n}\frac{1}{k^2}} \le \Pi_{k=0}^{n-1}\ell_{k,j}^{i_k} \le \frac{1}{2^n}e^{C_2\sum_{k=1}^{n}\frac{1}{k^2}} \le \frac{1}{2^n}e^{C_2 T}.
$$
That means
\begin{equation}
\label{eq:transverse measure on Omega_i}
\Pi_{k=0}^{n-1}\ell_{k,j}^{i_k} \asymp \frac{1}{2^n}.
\end{equation}

By (\ref{eq:transverse measure on Omega_i}), (\ref{eq:integral-over-Omega_1-union-Omega_2}) and $\ell (\alpha_n^j)\geq\frac{n^r}{2^n}$, the Dirichlet integral of $w$ over any level $n$ pair of pants $\Pi^{n,j}$ is less than
$$
K_1 [\max\{\mu_P^{n,j}, \mu_Q^{n,j}\}]^2 \cdot \frac{2^n}{n^r} \le K_2 \frac{1}{2^{2n}}\frac{2^n}{n^r} = K_2 \frac{1}{2^n \cdot n^r}
$$
for some $K_1,K_2>0$, where $\mu_P^{n,j}$ and $\mu_Q^{n,j}$ are the transverse measures of $\mathscr{F}$ on $p_n^j$ and $q_n^j$.

Since there are $2^n$ level $n$ pairs of pants, we have for some $K > 0$ that
$$
D_X(\mathscr{F})\leq K \sum_{n=1}^{\infty}\frac{1}{n^r}<\infty .
$$
By construction, both rays of every leaf in $\mathscr{F}$ leave every compact subset of $X_C$. Therefore $X_C\notin O_G$, and this finishes the proof of Theorem 1.1.

We proceed to assume that $\tilde{X}_C$ is the blooming Cantor tree surface with the properties stated in Theorem 1.2 . Assume, without loss of generality, no surface of genus at most $C$ and two boundaries is attached to $\alpha_0$. Refer to the front of the surface as $\tilde{X}$ for simplicity. Denote by $S_n^j$ the subset of $\tilde{X}$ that is the union of the front of a level $n$ pair of pants and the fronts of the sub-surfaces possibly attached to each of its level $n$ boundaries with genus at most $C$ and two geodesic boundaries for $n \ge 1$ and for each $1 \le j \le 2^{n}$.  Decompose $S_n^j$ into at most $8C+2$ sub-regions, half of which are examples of $P \cup R$ regions, while the rest are examples of $Q \cup S$ regions. From the estimates in (\ref{eq:integral-over-Omega_1-union-Omega_2}), Appendix B.1, and inequality $\ell_{n}^j \ge C_2 \frac{n^r}{2^{n}}$ for $n \ge 1$ and for each $1 \le j \le 2^{n+1}$,
$$
\iint_{S_n^j} (\frac{\partial v}{\partial x})^2 + (\frac{\partial v}{\partial y})^2 \lesssim \frac{2^{n}}{n^r}.
$$
Similar summation concludes that $\tilde{X}_C\notin O_G$.
$ \newline$
\vspace{-.25in}
\\
\begin{figure}[ht]
\centering
\vspace{-.5cm}
\includegraphics[width = 2.5in]{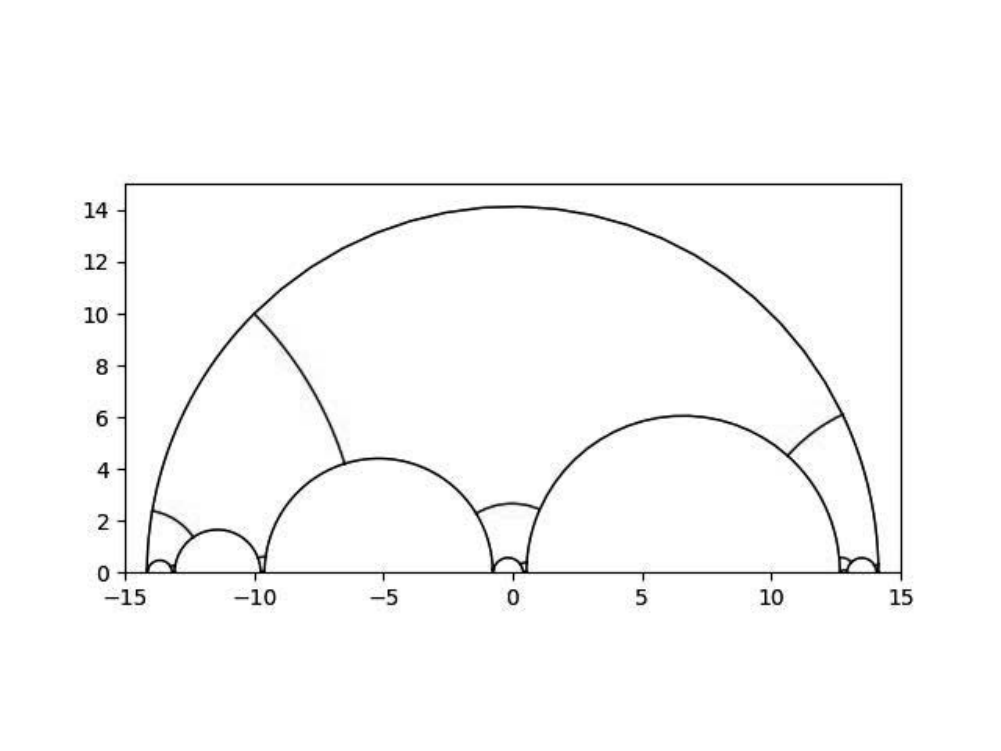}
\vspace{-1cm}
\caption{A computer-generated picture of an isomorphic lift to $\mathbb{H}$ of the front of a Cantor tree $X$. Lengths of cuffs vary between bounds in Theorem \ref{thm:Cantor_trees_non_parabolic}, where $r = 1.5$. The picture shows its geodesic flow escapes to infinity along every end due to the sizes of its cuffs. In the picture, $X$ is not parabolic since its Dirichlet integral is finite. There is a large amount of space on $X$ for the geodesics to escape towards many components of its end space $\delta_{\infty} X$.}
\label{fig:computer-generated picture of front of Cantor tree}
\end{figure}
$ \newline$
\vspace{-.45in}

\section*{Appendix A}
\indent Let $\Pi$ be a geodesic pair of pants with boundaries $\alpha_1$, $\alpha_2$, and $\alpha_3$. Let $o_{i,j}$ be the orthogeodesic arc between $\alpha_i$ and $\alpha_j$, for $i,j \in \{1,2,3\}$ such that $i \neq j$. The union $o_{1,2} \cup o_{1,3} \cup o_{2,3}$ separates $\Pi$ into front and back hexagons $H_1$ and $H_2$ with geodesic boundaries. Let $a_1$ be the orthogeodesic from $\alpha_1$ to $o_{2,3}$ that separates $H_1$ into two right-angled pentagons and divides $\alpha_1$ in $H_1$ from left to right into the sub-arcs $p$ and $q$ (see Figure \ref{fig:Horizontal Foliation through Omega one and two}). Call $P_p$ the pentagon containing $p$ (see Figure \ref{fig:Horizontal Foliation through Omega one and two}), and $P_q$ the pentagon containing $q$. 
The orthogeodesic $b_1$ from $a_1 \cap o_{2,3}$ to $o_{1,2}$ divides $P_p$ into quadrilaterals $P$ and $R$ adjacent to $\alpha_1$ and $\alpha_2$. Let $o_P$ be the sub-arc of $o_{1,2}$ in $P$ and let $o_R$ be the sub-arc of $o_{1,2}$ in $R$ (see Figure \ref{fig:Horizontal Foliation through Omega one and two}).

Recall that the {\it relative lengths} associated with $\alpha_1$ in the front of $\Pi$
is
$$
\ell^0 = \frac{\ell(p)}{\frac{\ell(\alpha_1)}{2}} \ and \ \ell^1 = \frac{\ell(q)}{\frac{\ell(\alpha_1)}{2}}. 
$$
Let $\ell = \ell^0$ or $\ell = \ell^1$.

\vskip .2 cm

\noindent {\bf Lemma A.2.}
\label{lem:bounds for relative cuff lengths}
{\it Assume $\frac{C_2}{(n+1)^2} \ge \max\{\ell(\alpha_2), \ell(\alpha_3)\}$ for a positive constant $C_2$ and for some $n \ge 1$. Then, any relative length $\ell$ associated with $\alpha_1$ in the front of $\Pi$ satisfies
$$
\frac{1}{2}e^{-\frac{C_2}{(n+1)^2}} \le \ell \le \frac{1}{2}e^{\frac{C_2}{(n+1)^2}}.
$$}

\begin{proof} We require the following hyperbolic trigonometric inequalities. Observe that
$$
\frac{d}{d x}(\frac{\sinh x}{A + \cosh x})|_{x=0} = \frac{1}{A + 1} \ and \ \frac{d}{d x}(\tanh (\frac{x}{2A}))|_{x = 0} = \frac{1}{2A}.
$$
It is also true that
$$
\frac{1}{A+1} > \frac{1}{2A}
$$
if and only if $A>1$. Let $A > 0$. For small $x > 0$, we get $A>1$ if and only if
\begin{equation}
\label{eq:trig-greater-than}
\frac{\sinh {x}}{A + \cosh {x}} > \tanh{(\frac{x}{2A})},
\end{equation}
with equality if and only if $A=1$. Using a formula for right-angled pentagons \cite[Theorem 7.18.1]{Beardon},
\begin{equation}
\label{eq:pentagon_formula}
\tanh{\ell(p)}\cosh{\ell(o_{1,2})}\tanh{\frac{\ell(\alpha_{2})}{2}} = 1.
\end{equation}
A formula for right-angled hexagons \cite[Theorem 7.19.2]{Beardon}, gives
\begin{equation}
\label{eq:hexagon_formula}
\cosh{\ell(o_{1,2})} = \frac{\cosh{\frac{\ell(\alpha_{3})}{2}}}{\sinh{\frac{\ell(\alpha_{1})}{2}}\sinh{\frac{\ell(\alpha_{2})}{2}}} + \coth{\frac{\ell(\alpha_{1})}{2}}\coth{\frac{\ell(\alpha_{2})}{2}}.
\end{equation}
Use (\ref{eq:trig-greater-than}), (\ref{eq:pentagon_formula}), and (\ref{eq:hexagon_formula}) to obtain for the relative length $\ell^0$ that
$$
\ell^0 = \frac{2}{\ell(\alpha_{1})}\tanh^{-1}{(\frac{\sinh{(\frac{\ell(\alpha_{1})}{2})}}{\frac{\cosh{(\frac{\ell(\alpha_{3})}{2})}}{\cosh{(\frac{\ell(\alpha_{2})}{2})}} + \cosh{(\frac{\ell(\alpha_{1})}{2})}})} > \frac{1}{2}[\frac{\cosh{(\min{\{\frac{\ell(\alpha_{2})}{2}, \frac{\ell(\alpha_{3})}{2}\}})}}{\cosh{(\max{\{\frac{\ell(\alpha_{2})}{2}, \frac{\ell(\alpha_{3})}{2}\}})}}] > \frac{e^{-\max{\{\ell(\alpha_{2}), \ell(\alpha_{3})\}}}}{2}.
$$
Fix a positive constant $C_2$ and fix some value of $n \ge 1$. From the assumption that $\frac{C_2}{(n+1)^2} \ge \max\{\ell(\alpha_2), \ell(\alpha_3)\}$, the above is greater than or equal to
$$
\frac{1}{2}e^{-\frac{C_2}{(n+1)^2}}.
$$
Similarly, we obtain the same lower bound for $\ell^1$ for $C_2>0$ and for $n \ge 1$. Use (\ref{eq:trig-greater-than}) to obtain
$$
\ell^0 = \frac{l(p)}{\frac{\ell(\alpha_{1})}{2}} < \frac{2}{\ell(\alpha_{1})} \tanh^{-1}{(\tanh{([\frac{\cosh{(\max{\{\frac{\ell(\alpha_{2})}{2}, \frac{\ell(\alpha_{3})}{2}\}})}}{\cosh{(\min{\{\frac{\ell(\alpha_{2})}{2}, \frac{\ell(\alpha_{3})}{2}\}})}}]\frac{\ell(\alpha_{1})}{4})})} < \frac{1}{2}e^{\frac{C_2}{(n+1)^2}}
$$
for $C_2 > 0$ and for $n \ge 1$.
Get the same upper bound for $\ell^1$ similarly for $C_2>0$ and for $n \ge 1$.
\end{proof}

\vskip .2 cm

\noindent {\bf Remark A.3.}
    {\it Lemma A.2 is still true when $\ell(\alpha_1)$ and $\ell(\alpha_3)$ are close and $\alpha_2$ is a puncture.}

\vskip .2 cm

We will require the following estimates. Refer to the beginning of Appendix A and Figure \ref{fig:image of Lambert quadrilateral to a subset of complex plane} for the definitions of $\Pi$ and $f(\tilde{P})$. 

\vskip .2 cm

\noindent {\bf Lemma A.4.}
\label{lem:proof of inequalities to estimate integral}
{\it We prove inequalities 1-3 in $f(\tilde{P})$ and inequality 4 in $\Pi$.
\begin{enumerate}
    \item Assume $\max\{\ell(\alpha_1), \ell(\alpha_2)\} \le B$ for some $B>0$. Then $\ell(b_1) \le D$ for some $D>0$.
    \item Assume $0 \le x \le \ell(o_P)$ and $\max\{\ell(\alpha_1), \ell(\alpha_2)\} \le B$ for some $B>0$. Then
$$
\frac{1}{1-[\tanh{\ell(p)\cosh{x}}]^2} \le c,
$$
for some $c>0$. 
    \item Assume $0 \le x \le \ell(o_P)$ and $\max\{\ell(\alpha_1), \ell(\alpha_2)\} \le B$ for some $B>0$. Then
$$
\tanh^{-1}{(\tanh{\ell(p)}\cosh{x})} \le c \tanh{\ell(p)}\cosh{x},
$$
for some $c>0$.
    \item If $\cosh\frac{\ell(\alpha_3)}{2} \le 2$, then $\ell(p) \gtrsim \ell(\alpha_1).$
\end{enumerate}}

\begin{proof}
We require the following limit.
\begin{equation}
\label{eq:limit}
\lim_{x \rightarrow 0^{+}}\sinh^{-1}{(\sinh{x} \cosh{(\sinh^{-1}{(\coth{x})})})} = \sinh^{-1}{1} \approx 0.88137.
\end{equation}
Let $s$ be the length of the geodesic sub-arc of $o_{2,3}$ from $a_1$ to $\alpha_2$. By a formula \cite[Theorem 7.18.1]{Beardon} for right-angled pentagons,
$$
\sinh{\ell(a_1)}\sinh{\ell(p)} = \cosh{\frac{\ell(\alpha_2)}{2}} \ and \ \sinh{s}\sinh{\frac{\ell(\alpha_2)}{2}} = \cosh{\ell(p)}.
$$
A formula for Lambert quadrilaterals \cite[Theorem 2.3.1(v)]{Buser} gives
$$
\sinh{\ell(b_1)} = \sinh{\ell(p)}\cosh{\ell(a_1)} \ and \ \sinh{\ell(b_1)} = \sinh{\frac{\ell(\alpha_2)}{2}}\cosh{s}.
$$
Let $A = \max\{\frac{\ell(\alpha_2)}{2}, \ell(p)\}$. If $\ell(p) \le \frac{\ell(\alpha_2)}{2}$, then since $\cosh x$ and $\cosh(\sinh^{-1}x)$ are increasing functions,
\begin{equation}
\label{eq:first_b_1_eqn}
\ell(b_1) = \sinh^{-1}(\sinh(\frac{\ell(\alpha_2)}{2})\cosh(\sinh^{-1}(\frac{\cosh(\ell(p))}{\sinh(\frac{\ell(\alpha_2)}{2})}))) \le \sinh^{-1}(\sinh(A)\cosh(\sinh^{-1}(\coth(A)))).
\end{equation}
Otherwise, from $\cosh x$ and $\cosh(\sinh^{-1}x)$ being increasing,
\begin{equation}
\label{eq:second_b_1_eqn}
\ell(b_1) = \sinh^{-1}(\sinh(\ell(p))\cosh(\sinh^{-1}(\frac{\cosh(\frac{\ell(\alpha_2)}{2})}{\sinh(\ell(p))}))) \le \sinh^{-1}(\sinh(A)\cosh(\sinh^{-1}(\coth(A)))).
\end{equation}
Therefore, by (\ref{eq:first_b_1_eqn}), (\ref{eq:second_b_1_eqn}), the definition of $A$, and (\ref{eq:limit}) we obtain
\begin{equation}
\label{eq:b_1_bound}
\ell(b_1) \le \sinh^{-1}{(\sinh{\frac{\ell(\alpha_1)}{2}}\cosh{(\sinh^{-1}{(\coth{\frac{\ell(\alpha_1)}{2}})})})} \le D,
\end{equation}
for some $D>0$. From two formulas \cite[Theorem 2.3.1(i) and (iii)]{Buser} for Lambert quadrilaterals and basic computations,
\begin{equation}
\label{eq:tanh_eqn}
\tanh \ell(b_1) = \tanh \ell(p) \cosh \ell(o_P).
\end{equation}
Consider $x$ such that $0 \le x \le \ell(o_P)$. Use $\cosh x$ and $\tanh^{-1}x$ are increasing, (\ref{eq:tanh_eqn}), and (\ref{eq:b_1_bound}) to obtain
$$\tanh^{-1}{(\tanh{\ell(p)}\cosh{x})} \le \tanh^{-1}{(\tanh{\ell(p)}\cosh{\ell(o_P)})} = \ell(b_1) \le D,
$$
$$
\frac{1}{1-[\tanh{\ell(p)\cosh{x}}]^2} \le \cosh^2{D}.
$$
Observe that
\begin{equation}
\label{eq:inequalities_with_E}
    \tanh{\ell(p)}\cosh{x} \le \tanh{D} \ and \ \tanh^{-1}{t} \le Et
\end{equation}
for some $E > 0$ and for $t \le \tanh{D}$, where $E$ depends on how small $\tanh D$ is. The larger $\tanh D$ is, the larger $E$ must be to make the second inequality in (\ref{eq:inequalities_with_E}) work. This implies
$$
\tanh^{-1}{(\tanh{\ell(p)}\cosh{x})} \lesssim \tanh{\ell(p)}\cosh{x}.
$$
Let $c = \max\{\cosh^2(D), E\}$ to attain inequalities (2) and (3) in Lemma A.4. From the hyperbolic trigonometric identity \cite[Theorem 7.18.1]{Beardon} in pentagon $P_1$ (see Figure \ref{fig:Horizontal Foliation through Omega one and two}), we obtain
\begin{equation}
\label{eq:l(o_{1,2})}
\cosh \ell(o_{1,2})=\frac{1}{\tanh \ell(p) \tanh \frac{\ell(\alpha_2)}{2}}.
\end{equation}
From formula (\ref{eq:l(o_{1,2})}), if $\ell(\alpha_3) \le \ell(\alpha_2)$, we get that
$$
\ell(p) = \tanh^{-1}{\frac{\sinh{\frac{\ell(\alpha_1)}{2}}}{\frac{\cosh{\frac{\ell(\alpha_3)}{2}}}{\cosh{\frac{\ell(\alpha_2)}{2}}} + \cosh{\frac{\ell(\alpha_1)}{2}}}} \ge \tanh^{-1}{\frac{\sinh{\frac{\ell(\alpha_1)}{2}}}{1 + \cosh{\frac{\ell(\alpha_1)}{2}}}} = \frac{\ell(\alpha_1)}{4}.
$$
If $\ell(\alpha_2) < \ell(\alpha_3)$, using the inequality (\ref{eq:trig-greater-than}) and the assumption $\cosh\frac{\ell(\alpha_3)}{2} \le 2$ we get
$$
\ell(p) > \tanh^{-1}({\tanh{([\frac{\cosh{\frac{\ell(\alpha_2)}{2}}}{\cosh{\frac{\ell(\alpha_3)}{2}}}]\frac{\ell(\alpha_1)}{4})}}) \gtrsim \ell(\alpha_1).
$$
\end{proof}
$ \newline$
\vspace{-.65in}
\section*{Appendix B}
The sub-arc $o_P$ of the geodesic arc between two boundary cuffs of a pair of pants of a decomposition of $X$ in quadrilateral $P$ lifts to the arc $[i, e^{\ell(o_P)}i]$ on the y-axis. Each point $w$ in the lift $\tilde{P}$ of quadrilateral $P$ belongs to the lift $\tilde{\gamma}_{w_0}$ of some geodesic arc $\gamma_{w_0}$ orthogonal to $o_P$ with foot $w_0$ on $o_P$. Note that $\tilde{\gamma}_{w_0}$ is a geodesic arc orthogonal to the y-axis and to its right. Map $\tilde{\gamma}_{w_0}$ by an isometry $f$ to the Euclidean segment orthogonal to $[0, \ell(o_P)]$ and below the x-axis for each $w_0$ in $o_P$ (see Figure \ref{fig:image of Lambert quadrilateral to a subset of complex plane}). That completely defines $f$ on $\tilde{P}$. We define $f$ similarly on $\tilde{Q}$, $\tilde{R}$, and $\tilde{S}$.

Consider the inverse $g = f^{-1}$ of the diffeomorphism $f$ defined on the lift $\tilde{P}$, $\tilde{Q}$, $\tilde{R}$, or $\tilde{S}$ of $P$, $Q$, $R$, or $S$ in the front of a pair of pants $\Pi$. In this section we explain how $g$ is quasiconformal with a quasiconformal constant that is bounded above by $k_0 = \frac{1 + \csch^2 (D)}{\coth (D) \csch (D)}$.

\vskip .2 cm

\noindent {\bf Proposition B.1.}
\label{lem:quasiconformal constant}
{\it The diffeomorphisms $f$ from $\tilde{P}$, $\tilde{Q}$, $\tilde{R}$, and $\tilde{S}$ to $\mathbb{C}$ is quasiconformal with quasiconformal constant bounded above by $k_0 = \frac{1 + \csch^2 (D)}{\coth (D) \csch (D)}$.}

\begin{proof} We obtain a diffeomorphism $f$ from $\tilde{P}$ to $\mathbb{C}$ with quasiconformal constant bounded by $k_0$. Extend the result generally to every additional case. The inverse map $g$ of $f$ sends a point $z = x + iy$ to $w$ such that $|w| = e^x$ and $\arg(w) = \tan^{-1}(\csch (-y))$. Thus, $g$ is defined at $z = x + iy$ by
$$
g(z) = f^{-1}(z) = e^xe^{i\tan^{-1}(\csch (-y))} = e^{\frac{z + \bar{z}}{2}}e^{i\tan^{-1}(\csch(\frac{\bar{z} - z}{2i}))}.
$$
The dilatation of $g$ at $z = x + iy$, denoted by $K(z)$, is
$$
K(z) = \frac{1 + \csch^2 (-y)}{\coth (-y) \csch (-y)}.
$$
The value of $-y$ in $f(\tilde{P})$ is bounded above by $\ell(b_1)$, which is bounded above by a positive constant $D$ (see Lemma A.4(1)). Note $\frac{1 + \csch^2 (-y)}{\coth (-y) \csch (-y)}$ is decreasing with respect to negative $y$. Thus, $g = f^{-1}$ is quasiconformal with quasiconformal constant bounded above by $k_0 = \frac{1 + \csch^2 (D)}{\coth (D) \csch (D)}$.
\end{proof}


\end{document}